\documentclass[12pt]{article}
\usepackage{mathrsfs}
\usepackage{amscd}
\usepackage{amsmath,amsfonts,amssymb,amscd}
\usepackage{indentfirst,graphics,epsfig,psfrag}
\input{epsf}
\usepackage{ifpdf}
\usepackage{enumerate}
\usepackage{appendix}
\usepackage{enumerate}

\usepackage{lineno}

\makeatletter \@addtoreset{figure}{section} \makeatother
\makeatletter
\long\def\@makecaption#1#2{%
   \vskip 10\p@
   \setbox\@tempboxa\hbox{{#1}\ \ #2}%
   \ifdim \wd\@tempboxa >\hsize
       {#1}\ \ #2\par
   \else
       \hbox to\hsize{\hfil\box\@tempboxa\hfil}%
   \fi}
\makeatother

\newtheorem{thm}{Theorem}[section]
\newtheorem{cor}[thm]{Corollary}

\newcommand{\qed}{{\hfill\rule{3pt}{7pt}}}

\def\qed{\hfill \rule{4pt}{7pt}}

\begin{document}
\title{\bf The rainbow $k$-connectivity\\
of two classes of graphs\footnote {Supported by NSFC, PCSIRT and the
``973" program. }}
\author{
\small  Xueliang Li,  Yuefang Sun\\
\small Center for Combinatorics and LPMC-TJKLC\\
\small Nankai University, Tianjin 300071, P.R. China\\
\small E-mail: lxl@nankai.edu.cn; syf@cfc.nankai.edu.cn
  }
\date{}
\maketitle
\begin{abstract}
A path in an edge-colored graph $G$, where adjacent edges may be
colored the same, is called a rainbow path if no two edges of $G$
are colored the same. For a $\kappa$-connected graph $G$ and an
integer $k$ with $1\leq k\leq \kappa$, the rainbow $k$-connectivity
$rc_k(G)$ of $G$ is defined as the minimum integer $j$ for which
there exists a $j$-edge-coloring of $G$ such that every two distinct
vertices of $G$ are connected by $k$ internally disjoint rainbow
paths. Let $G$ be a complete $(\ell+1)$-partite graph with $\ell$
parts of size $r$ and one part of size $p$ where $0\leq p <r$ (in
the case $p=0$, $G$ is a complete $\ell$-partite graph with each
part of size $r$). This paper is to investigate the rainbow
$k$-connectivity of $G$. We show that for every pair of integers
$k\geq 2$ and $r\geq 1$, there is an integer $f(k,r)$ such that if
$\ell\geq f(k,r)$, then $rc_k(G)=2$. As a consequence, we improve
the upper bound of $f(k)$ from $(k+1)^2$ to $ck^{\frac{3}{2}}+C$,
where $0<c<1$, $C=o(k^{\frac{3}{2}})$, and $f(k)$ is
the integer such that if $n \geq f(k)$ then $rc_k(K_n)=2$.\\[2mm]
\end{abstract}

\section{Introduction}

All graphs considered in this paper are simple, finite and
undirected. Let $G$ be a nontrivial connected graph with an edge
coloring $c: E(G)\rightarrow \{1,2,\cdots,k\}$, $k\in \mathbb{N}$,
where adjacent edges may be colored the same. A path of $G$ is
called $rainbow$ if no two edges of it are colored the same. An edge
colored graph $G$ is said to be $rainbow ~connected$ if for any two
vertices of $G$ there is a rainbow path of $G$ connecting them.
Clearly, if a graph is rainbow connected, it must be connected.
Conversely, any connected graph has a trivial edge coloring that
makes it rainbow connected, i.e., the coloring such that each edge
has a distinct color. Thus, Chartrand et al. \cite{Chartrand 2}
defined the $rainbow~connection~number$ of a connected graph $G$,
denoted by $rc(G)$, as the smallest number of colors by using which
there is an edge coloring of $G$ such that $G$ is rainbow connected.

Suppose that $G$ is a $\kappa$-connected graph with $\kappa \geq 1$.
It then follows from a well-known theorem of Whitney that for every
integer $k$ with $1\leq k\leq \kappa$ and every two distinct
vertices $u$ and $v$ of $G$, the graph $G$ contains $k$ internally
disjoint $u-v$ paths. Chartrand et al. \cite{Chartrand 2} defined
the $rainbow$~ $k$-$connectivity$ $rc_k(G)$ of $G$ to be the minimum
integer $j$ for which there exists a $j$-edge-coloring of $G$ such
that for every two distinct vertices $u$ and $v$ of $G$, there exist
at least $k$ internally disjoint $u-v$ rainbow paths. Thus by the
definition, $rc_1(G)=rc(G)$ is the rainbow connection number of $G$.
By coloring the edges of $G$ with distinct colors, we see that every
two vertices of $G$ are connected by $\kappa$ internally disjoint
rainbow paths and that $rc_k(G)$ is defined for every integer $k$
with $1 \leq k \leq \kappa$. Furthermore, $rc_k(G)\leq rc_j(G)$ for
$1\leq k\leq j\leq \kappa$.

The concept of rainbow $k$-connectivity has applications in
communication networks. For details we refer to \cite{Chartrand 2}
and \cite{Ericksen}.

In \cite{Chartrand 2},  Chartrand et al. studied the rainbow
$k$-connectivity of the complete graph $K_n$ for various pairs $k,
n$ of integers. It was shown in \cite{Chartrand 2} that for every
integer $k\geq 2$, there exists an integer $f(k)$ such that
$rc_k(K_n)=2$ for every integer $n\geq f(k)$. They also investigated
the rainbow $k$-connectivity of $r$-regular complete bipartite
graphs for some pairs $k, r$ of integers with $2\leq k\leq r$, and
they showed that for every integer $k\geq 2$, there exists an
integer $r$ such that $rc_k(K_{r,r})=3$.

This paper is to continue their investigation. We improve the upper
bound of $f(k)$ from $(k+1)^2$ to $ck^{\frac{3}{2}}+C$ (here $0<c<1$
and $C=o(k^{\frac{3}{2}})$), i.e., from $O(k^2)$ to
$O(k^{\frac{3}{2}})$, a considerable improvement. For notations and
terminology not defined here, we refer to \cite{Bondy}.

\section{Main Results}

In \cite{Chartrand 2} the following result was obtained.

\begin{thm}\label{thm1}(\cite{Chartrand 2}) For every integer $k\geq
2$, there exists an integer $f(k)$ such that if $n\geq f(k)$, then
$rc_k(K_n)=2$. \qed
\end{thm}

The authors in \cite{Chartrand 2} then obtained an upper bound
$(k+1)^2$ for $f(k)$, namely $f(k)\leq (k+1)^2$. In the following we
will give a better upper bound for $f(k)$.

Let $K_{\ell[r]}$ denote the complete $\ell$-partite graph each part
of which contains $r$ elements. Let $V=\bigcup_{s=1}^\ell{V_s}$
where $V_s=\{u_{s,1},u_{s,2},\cdots,u_{s,r} \}$ $(1\leq s\leq \ell)$
is the vertex set of each part, and
$U_j=\{u_{1,j},u_{2,j},\cdots,u_{l,j}\}$ $(1\leq j\leq r)$.

We consider the graph $G=K_{\ell^2[r]}$ where $\ell$ is an integer,
and give $G$ a 2-edge-coloring as follows: Let
$\overline{V_i}=\bigcup_{t=1}^{\ell}{V_{(i-1)\ell+t}}$ where $1\leq
i\leq \ell$. We know that the subgraph, denoted by $G_j$, of $G$
induced by $U_j$ is the complete graph $K_{\ell^2}$. Let
$U_j=U_{1,j} \cup \cdots \cup U_{\ell,j}$ where
$U_{i,j}=\{u_{(i-1)\ell+1,j},\cdots,u_{i\ell,j} \}$. Then the
subgraph, denoted by $G_{i,j}$, of $G$ induced by $U_{i,j}$ is the
complete graph $K_\ell$. Similar to the edge coloring in the proof
of Theorem \ref{thm1} in \cite{Chartrand 2}, we first give a
2-edge-coloring of each $G_j$ $(1\leq j\leq r)$ as follows: We
assign the edge $uv$ of $G$ the color 1 if either $uv \in
E(G_{i,j})$ for some $i(1\leq i\leq \ell)$ or if
$uv=u_{(i_1-1)\ell+t,j}~u_{(i_2-1)\ell+t,j}$ for some $i_1,i_2,t$
with $1\leq i_1,i_2,t\leq \ell$ and $i_1\neq i_2$. All other edges
of $G_j$ are assigned the color 2. With a similar argument to
Theorem \ref{thm1} in \cite{Chartrand 2}, there are $\ell-1$
disjoint rainbow $u-v$ paths for any two vertices $u, v$ in the
subgraph $G_j$ including 1 path of length 1 and $\ell-2$ paths of
length 2.

For other edges, that is, the edges between distinct $G_j$s, we use
above two colors as follows: Let $uv=u_{i_1,j_1}u_{i_2,j_2}$ where
$u_{i_1,j_1} \in G_{j_1}, u_{i_2,j_2}\in G_{j_2}$, we give it the
color different from the color of edge $u_{i_1,j_1}u_{i_2,j_1}$ of
graph $G_{j_1}$ or $u_{i_1,j_2}u_{i_2,j_2}$ of graph $G_{j_2}$
(since edges $u_{i_1,j_1}u_{i_2,j_1}$ and $u_{i_1,j_2}u_{i_2,j_2}$
have the same color by above coloring).

Next we will count the number of disjoint rainbow $u-v$ paths for
any two vertices $u,v$ of $G$ with above edge coloring. Without loss
of generality, let $u=u_{1,1}$. We consider four cases:

\textbf{Case 1}. $v=u_{1,j}$ $(2 \leq j\leq r)$, that is, $u$ and
$v$ are in the same part.

By above coloring, edges $u_{1,1}u_{s,j}$ and $u_{s,j}u_{1,j}$ have
distinct colors where $2 \leq s\leq \ell^2$, and so path
$u_{1,1},u_{s,j},u_{1,j}$ is a rainbow $u-v$ path. Similarly, path
$u_{1,1},u_{s,1},u_{1,j}$ is a rainbow $u-v$ path for each $2 \leq
s\leq \ell^2$. And any two such rainbow paths are disjoint, and so
there are at least $2(\ell^2-1)$ disjoint rainbow $u-v$ paths in
$G$.

\textbf{Case 2}. $v=u_{s,1}$ $(2\leq s\leq \ell^2)$, that is, $u$
and $v$ are in the same $G_j$ (here $j=1$).

Then in subgraph $G_1$ there are $\ell-1$ disjoint rainbow $u-v$
paths, and in each subgraph $G_j$ where $(2\leq j\leq r)$, there are
$\ell-2$ disjoint rainbow $u_{1,j}-u_{s,j}$ paths of length 2
$u_{1,j},y,u_{s,j}$ for some $y$s. Since the color of $u_{1,j}y$,
$yu_{s,j}$ is different from that of $u_{1,1}y$, $yu_{s,1}$,
respectively, paths $u_{1,1},y,u_{s,1}$ are $\ell-2$ disjoint
rainbow paths in each $G_j$ where $(2\leq j\leq r)$. So there are
totally $r(\ell-2)+1$ disjoint rainbow $u-v$ paths in $G$.

\textbf{Case 3}. $v=u_{s,j_0}$ where $2\leq s\leq \ell, 2\leq
j_0\leq r$, that is, $u$ and $v$ belong to the same $\overline{V_i}$
(here $i=1$) but in distinct parts and subgraphs $G_j$s.

At first, the edge $u_{1,1}u_{s,j_0}$ is a rainbow $u-v$ path.

Next we consider the monochromatic $u_{1,1}-u_{s,1}$ path
$u_{1,1},y, u_{s,1}$. Then each path $u_{1,1},y, u_{s,j_0}$ is a
rainbow $u_{1,1}-u_{s,j_0}$ path of length 2 where $2\leq s\leq
\ell, 2\leq j_0\leq r$. In $U_1$, we choose $y=u_{(i-1)\ell+t,1}$
where $1\leq i\leq \ell, 2\leq t\leq \ell$ and $t\neq s$. So there
are $\ell(\ell-2)$ such paths.

Similarly, there are $\ell(\ell-2)$ disjoint monochromatic
$u_{1,j_0}-u_{s,j_0}$ paths $u_{1,j_0},y,u_{s,j_0}$ where $y\in
U_{j_0}$ and we can obtain other $\ell(\ell-2)$ disjoint rainbow
$u-v$ paths.

Moreover, in each subgraph $G_j$ where $2\leq j\leq r$ and $j\neq
j_0$, there are $\ell-2$ disjoint rainbow $u_{1,j}-u_{s,j}$ path
$u_{1,j},y,u_{s,j}$ of length 2. Since the color of $u_{1,j}y$ and
$yu_{s,j}$ is different from that of $u_{1,1}y$ and $yu_{s,j_0}$,
respectively, paths $u_{1,1},y,u_{s,j_0}$ are disjoint rainbow $u-v$
paths.

So the number of disjoint rainbow $u-v$ paths in $G$ is
$1+2\ell(\ell-2)+(r-2)(\ell-2)=1+(2\ell+r-2)(\ell-2)$.

\textbf{Case 4}. $v=u_{{(i_0-1)\ell+t_0},j_0}$ where $2\leq i_0\leq
\ell, 1\leq t_0\leq \ell, 2\leq j_0\leq r$, this is, $v$ and $u$ are
in distinct $\overline{V_i}$s (here $v$ is not in $\overline{V_1}$).

\textbf{Subcase 4.1}. $t_0=1$.

Then edge $uv$ is a rainbow path. In subgraph $G_1$, we find the
monochromatic $u-u_{(i_0-1)\ell+1,1}$ path
$u,y,u_{(i_0-1)\ell+1,1}$. We choose any vertex of $U_1\setminus
\{U_{1,1}\cup U_{i,1}\}$ as $y$. Since the color of edge
$yu_{(i_0-1)\ell+1,1}$ is different from that of $yv$, each path
$u,y,v$ is rainbow, and there are $\ell^2-2\ell$ disjoint rainbow
$u-v$ paths.

Similarly, in subgraph $G_{j_0}$, we can find monochromatic paths
$u_{1,j_0},y,u_{(i_0-1)\ell+1,j_0}$, and get other $\ell^2-2\ell$
disjoint rainbow $u-v$ paths.

In each subgraph $G_j$ where $j\neq 1,j_0$, we find the disjoint
rainbow $u_{1,j}-u_{(i_0-1)\ell+1,j}$ paths
$u_{1,j},y,u_{(i_0-1)\ell+1,j}$ of length 2. Since the color of
$u_{1,j}y$ and $yu_{(i_0-1)\ell+1,j}$ is different from that of
$u_{1,1}y$ and $yu_{(i_0-1)\ell+1,1}$, respectively, there are
$\ell-2$ disjoint rainbow $u-v$ paths
$u_{1,1},y,u_{(i_0-1)\ell+1,1}$.

So the number of disjoint rainbow $u-v$ paths in $G$ is
$1+2\ell(\ell-2)+(r-2)(\ell-2)=1+(2\ell+r-2)(\ell-2)$.

\textbf{Subcase 4.2}. $t_0\neq 1$.

Then, edge $uv$ is a rainbow path. In subgraph $G_1$, we find the
monochromatic $u-u_{(i_0-1)\ell+t_0,1}$ path
$u,y,u_{(i_0-1)\ell+t_0,1}$. Let $y=u_{(i-1)\ell+t,1}$ where $i\neq
1, i_0, t\neq 1,t_0$ or $y=u_{t_0,1}$ or $u_{(i_0-1)\ell+1,1}$.
Since the color of edge $yu_{(i_0-1)\ell+t_0,1}$ is different from
that of $yv$, each path $u,y,v$ is rainbow, and so there are
$(\ell-2)(\ell-2)$+2 disjoint rainbow $u-v$ paths.

Similarly, in subgraph $G_{j_0}$, we can find monochromatic paths
$u_{1,j_0},y,u_{(i_0-1)\ell+1,j_0}$, and get other
$(\ell-2)(\ell-2)+2$ disjoint rainbow $u-v$ paths.

In each subgraph $G_j$ where $j\neq 1,j_0$, we find the disjoint
rainbow $u_{1,j}-u_{(i_0-1)\ell+t_0,j}$ paths
$u_{1,j},y,u_{(i_0-1)\ell+t_0,j}$ of length 2. Since the color of
$u_{1,j}y$ and $yu_{(i_0-1)l+t_0,j}$ is different from that of
$u_{1,1}y$ and $yu_{(i_0-1)l+t_0,j_0}$, respectively, there are
$\ell-2$ disjoint rainbow $u-v$ paths
$u_{1,1},y,u_{(i_0-1)\ell+t_0,j_0}$ (here $y\in G_j$ where $j\neq
1,j_0$).

So the number of disjoint rainbow $u-v$ paths in $G$ is
$5+2(\ell-2)(\ell-2)+(r-2)(\ell-2)=5+(2\ell+r-6)(\ell-2)$.

We know here $\ell\geq  2$, and so $5+(2\ell+r-6)(\ell-2)
>r(\ell-2)+1$. Clearly, $1+(2\ell+r-2)(\ell-2) >r(\ell-2)+1$. Let $\ell_0=\lceil
\max \{\sqrt{\frac{k}{2}+1},\frac{k-1}{r}+2 \}\rceil$. We then have
$2(\ell_0^2-1)\geq k$, $r(\ell_0-2)+1\geq k$. So, our following
theorem holds.

\begin{thm}\label{thm2}
For every integer $k\geq 2$, there exists an integer $\ell_0=\lceil
\max \{\sqrt{\frac{k}{2}+1},\frac{k-1}{r}+2 \}\rceil$, such that if
$\ell\geq \ell_0$, then $rc_k(K_{\ell^2[r]})=2.$ \qed
\end{thm}

\begin{figure}[!hbpt]
\begin{center}
\includegraphics[scale=1.000000]{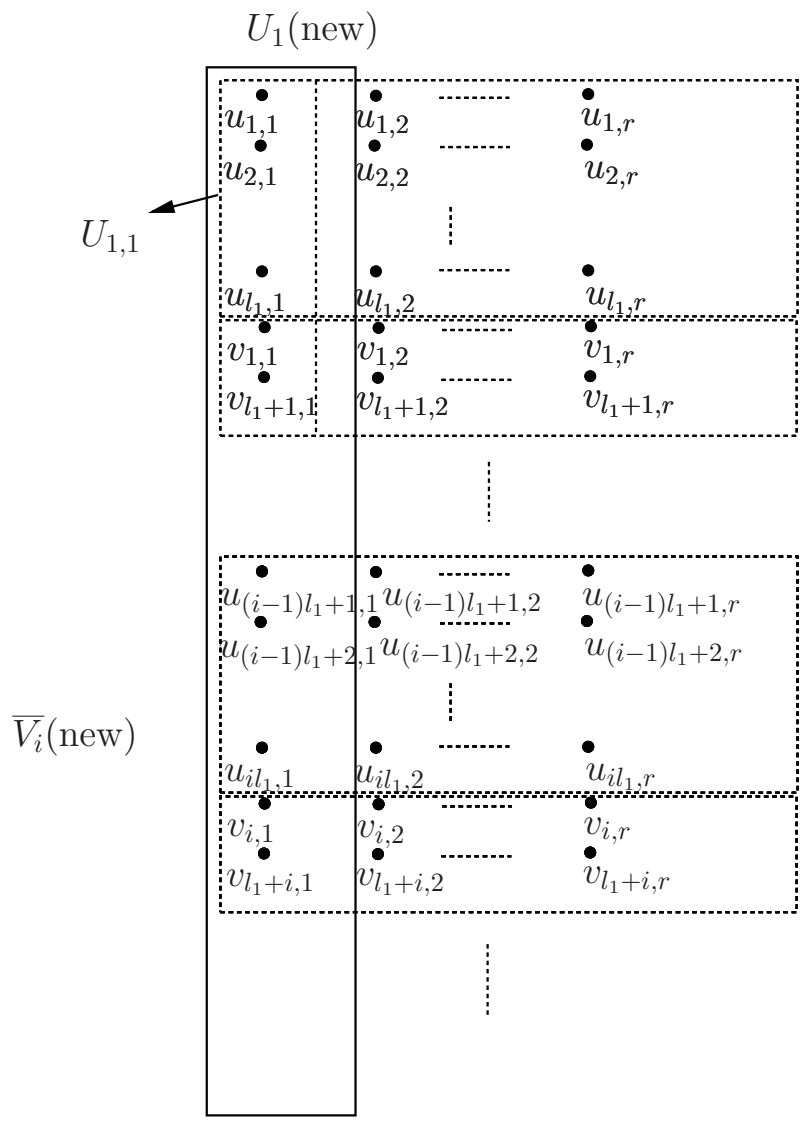}
\end{center}
\caption{The figure of Theorem \ref{thm3}.}\label{figure1}
\end{figure}

Next we will consider the case for general complete multipartite
graph $G\cong K_{\ell[r]}$ with equal part, where $\ell\geq
\ell_0^2$. Let $\ell_1$ be the integer satisfying $\ell_1^2 \leq
\ell\leq (\ell_1+1)^2$ where $\ell_1 \geq \ell_0$. Then by our
Theorem \ref{thm2}, we need to consider the case $1\leq
q=\ell-\ell_0^2 \leq 2\ell_1$. Now $G$ can be obtained from
$K_{{\ell_1^2}[r]}$ by adding $q$ new parts $P_i$~($1\leq i\leq
q$). Let the vertex set of each corresponding new part be
$V(P_i)=\{v_{i,1},v_{i,2},\cdots,v_{i,r} \}~(1\leq i\leq q)$. We
now update $\overline{V_i}$ as follows: If $1\leq q\leq \ell_1$,
let the new $\overline{V_i}$ be the union of the old
$\overline{V_i}$ and $V(P_i)$; If $\ell_1< q \leq 2{\ell_1}$, let
the new $\overline{V_i}$ be the union of the old $\overline{V_i}$
and $V(P_i)\cup V(P_{\ell_1+i})$ (if there exists) as shown in the
Figure \ref{figure1}. Similarly, we update $U_j~(1\leq j\leq r)$
(also $G_j, U_{i,j}, G_{i,j}$) as shown in the Figure
\ref{figure1} (for example, the new $U_1$).

Similar to the proof of Theorem \ref{thm1} in \cite{Chartrand 2}, we
give the 2-edge-coloring of graph $G$ as follows: For each
$G_j~(1\leq j\leq r)$, we assign the edge $uv$ of $G$ the color 1 if
either $uv \in E(G_{i,j})$ for some $i~(1\leq i\leq \ell_1)$ or if
$uv=u_{(i_1-1)\ell_1+t,j}~u_{(i_2-1)\ell_1+t,j}$ or
$uv=v_{i_1,j}~v_{i_2,j}$ or $uv=v_{\ell_1+i_1,j}~v_{\ell_1+i_2,j}$
for some $i_1,i_2,t$ with $1\leq i_1,i_2,t\leq \ell_1$ and $i_1\neq
i_2$. All other edges of $G_j$ are assigned the color 2.

For other edges, that is, edges between distinct $G_j$s, we use
above two colors as follows: If $uv=u_{i_1,j_1}~u_{i_2,j_2}$ where
$u_{i_1,j_1} \in G_{j_1}, ~u_{i_2,j_2}\in G_{j_2}$, we give it the
color different from that of edge $u_{i_1,j_1}~u_{i_2,j_1}$ of graph
$G_{j_1}$ or $u_{i_1,j_2}~u_{i_2,j_2}$ of graph $G_{j_2}$ (since
edges $u_{i_1,j_1}~u_{i_2,j_1}$ and $u_{i_1,j_2}~u_{i_2,j_2}$ have
the same color by above coloring); Similarly, color the edge for the
cases that $uv=u_{i_1,j_1}v_{i_2,j_2}$ and
$uv=v_{i_1,j_1}v_{i_2,j_2}$.

By the proof of our Theorem \ref{thm2}, there are at least $k$
disjoint rainbow paths connecting any two vertices $u,v \in
V(G)\setminus {\bigcup_{i=1}^q {V(P_i)}}$. We need to consider the
paths between $u,v$ where $u\in {\bigcup_{i=1}^q {V(P_i)}}$ and
$v\in V(G)$. Without loss of generality, let $u=v_{1,1}$, and we
consider four cases, which are similar to the four cases in the
proof of our Theorem \ref{thm2}.

\textbf{Case 1}.  $v\in V(P_1)$, that is, $u$ and $v=v_{1,j_0}$
$2\leq j_0\leq r$ are in the same part (here the part is $P_1$).

$G$ contains the disjoint rainbow $u-v$ paths $v_{1,1},y,v_{1,j_0}$,
and we let $y$ be the elements of $U_1 \cup U_{j_0} \setminus
\{v_{1,1}, v_{1,j_0} \}$. So the number of these paths is at least
$2\ell_1^2 \geq 2\ell_0^2 >2(\ell_0^2-1)\geq k$.

\textbf{Case 2}. $v\in U_1$, that is, $v$ and $u$ are in the same
subgraph $G_1$,

In $G_1$, there are at least $\ell_1-1$ rainbow $u-v$ paths; In each
$G_j$ $(2 \leq j\leq r)$ there are at least $\ell_1-2$ rainbow
$u'-v'$ paths of length 2: $u', y, v'$ where $u'$, $v'$ is in the
same part as $u$, $v$, respectively, and we get at least $\ell_1-2$
rainbow $u-v$ paths: $u, y, v$, where $y \in G_j$ $(2 \leq j\leq
r)$.

So the number of disjoint rainbow $u-v$ paths is at least
$1+r(\ell_1-2)\geq 1+r(\ell_0-2)\geq k$.

\textbf{Case 3}. $v \in \overline{V_1}$ but not in $U_1 \cup
V(P_1)$. Without loss of generality, let $v \in U_{j_0}$.

Then, edge $uv$ is a $u-v$ path.

\textbf{Subcase 3.1}. $v=v_{l+1,j_0}$.

Then similar to Case 3 in the proof of of Theorem \ref{thm2}, we
find the monochromatic $u-v_{l+1,1}$ path of length 2 in subgraph
$G_1$ and get at least $\ell_1^2$ disjoint rainbow $u-v$ paths.
Similarly, there are other $\ell_1^2$ disjoint rainbow $u-v$ paths
by finding the monochromatic $v_{1,j_0}-v_{\ell+1,j_0}$ paths in
$G_{j_0}$.

In each $G_j$ $(j\neq 1,j_0)$, similar to Case 3 in the proof of
Theorem \ref{thm2}, we get at least $\ell_1-2$ disjoint rainbow
$u-v$ paths of length 2 by finding rainbow $v_{1,j}-v_{\ell+1,j}$
paths in $G_j$.

So the number of disjoint rainbow $u-v$ paths is at least
$1+2\ell_1^2+
(r-2)(\ell_1-2)>1+2\ell_1(\ell_1-2)+(r-2)(\ell_1-2)=1+(2\ell_1+r-2)(\ell_1-2)\geq
1+(2\ell_0+r-2)(\ell_0-2)\geq k$.

\textbf{Subcase 3.2}. $v=u_{s_0,j_0}$ where $1\leq s_0 \leq \ell_1,
2\leq j_0\leq r$.

With a similar procedure to Subcase 3.1, we can get at least
$1+2\ell_1(\ell_1-1)+ (r-2)(\ell_1-2)> 1+(2\ell_1+r-2)(\ell_1-2)\geq
1+(2\ell_0+r-2)(\ell_0-2)\geq k$ rainbow disjoint $u-v$ paths in
graph $G$.

\textbf{Case 4}. $v$ is not in $\overline{V_1}$ and $U_1$, without
loss of generality, let $v \in \overline{V_2}\setminus U_1$. We
consider two subcases similar to Case 4 in the proof of Theorem
\ref{thm2}.

\textbf{Case 4.1}. $v=v_{2,j_0}$ where $2\leq j_0\leq r$.

Similar to Subcase 4.1 in the proof of Theorem \ref{thm2}, we can
get at least
$1+2\ell_1(\ell_1-2)+(r-2)(\ell_1-2)=1+(2\ell_1+r-2)(\ell_1-2) \geq
1+(2\ell_0+r-2)(\ell_0-2)\geq k$ disjoint rainbow $u-v$ paths in
$G$.

\textbf{Case 4.2}. $v\neq v_{2,j_0}$.

Similar to Subcase 4.2 in the proof of Theorem \ref{thm2}, we can
get at least
$1+2(\ell_1-1)(\ell_1-2)+(r-2)(\ell_1-2)=1+(2\ell_1+r-4)(\ell_1-2)
\geq 1+r(\ell_1-2)\geq 1+r(\ell_0-2)\geq k$ (since $\ell_1\geq 2$)
disjoint rainbow $u-v$ paths in $G$.

With above discussion, we have our following theorem:
\begin{thm}\label{thm3}
For every integer $k\geq 2$, there exists an integer $\ell_0=\lceil
\max \{\sqrt{\frac{k}{2}+1},\frac{k-1}{r}+2 \}\rceil$, such that if
$\ell\geq \ell_0^2$, then $rc_k(K_{\ell[r]})=2.$ \qed
\end{thm}

We now obtain a new graph $G'$ from $G$ by adding $p$ new vertices:
$w_1,w_2,\cdots,w_p$ $(1 \leq p\leq r-1)$ and edge $vw_j$ where
$1\leq j\leq p$ and $v \in V(G)$. Color the new edges as follows:
Give color 2 to edges between $U_j$ and $w_j$, and color 1 to other
edges.

Now we count the number of rainbow disjoint $u-v$ paths between any
two vertices of $G$. We need to consider the case that $u=w_j, v\in
V(G')$ $(1\leq j\leq p)$. Without loss of generality, let $u=w_1$.

\textbf{Case 1}. $v \in V(G)$.

\textbf{Subcase 1.1}. $v \in U_1$, without loss of generality, let
$v=u_{1,1}$.

Then $G$ contains the $u-v$ path $u,v$ as well as the $u-v$ rainbow
paths $w_1,y,u_{1,1}$ where $y \in {\overline{V_1}\setminus
\{u_{1,1},\cdots,u_{1,r}\}}$, or $y=u_{(i-1){\ell_1}+1,j} (2\leq
i\leq \ell_1, 1\leq j\leq r )$. So the number of disjoint rainbow
$u-v$ paths is at least $1+2r(\ell_1-1)
>1+r(\ell_0-2) \geq k$.

\textbf{Subcase 1.2}. $v$ is not in $U_1$, without loss of
generality, let $v=u_{1,2}$.

Then $G$ contains the $u-v$ path $u,v$ as well as the $u-v$ rainbow
paths $w_1,y,u_{1,2}$ where $y \in (U_1 \cup U_2) \setminus (U_{1,1}
\cup U_{1,2}\cup \{u_{(i-1){\ell_1+1,j_0}}\}_{i=1}^{\ell_1})$ where
$j_0=1,2$, or $y \in \{u_{(i-1){i_1}+1,j_1} \}_{i=2}^{\ell_1}$ where
$j_1\neq 1,2$, or $y= u_{s,j}$ where $2\leq s\leq \ell_1, j\neq 1,
2$. So the number of disjoint rainbow $u-v$ paths is at least
$1+2({\ell_1}-1)(\ell_1-1)+2(\ell_1-1)(r-2)=1+(\ell_1-1)(2\ell_1+2r-6)
>1+(\ell_1-2)(2\ell_1+2r-6)\geq 1+(\ell_1-2)r \geq k$(here we let $r\geq 2$).

\textbf{Case 2}. $v=w_j$ where $j\neq 1$.

Then $G$ contains the $u-v$ rainbow paths $w_1,y,w_j$ where $y\in
U_1 \cup U_j$. So the number of disjoint rainbow $u-v$ paths is at
least $2\ell_1^2 \geq 2\ell_0^2 >2(\ell_0^2-1) \geq k$.

If $r=1$, then $G'=G$ is complete, so by Theorem \ref{thm1} of
\cite{Chartrand 2} and above discussion, we have the following
theorem:
\begin{thm}\label{thm3}
Let $G$ be a complete $(\ell+1)$-partite graph with $\ell$ parts of
size $r \geq 1$ and a part of size $p$ where $1 \leq p <r$. Then for
every integer $k\geq 2$, if $\ell\geq \ell_0^2$ where $\ell_0=\lceil
\max \{\sqrt{\frac{k}{2}+1},\frac{k-1}{r}+2 \}\rceil$, then
$rc_k(G)=2.$ \qed
\end{thm}

We know that the graph in Theorem \ref{thm3} is a spanning
subgraph of the complete $K_{n}$ with $n=\ell r+p$ and $\ell\geq
\ell_0^2$. Now we let $r=r_0=\lceil \sqrt{2k}\rceil$, and $G$ be
the complete graph $K_n$ with $n \geq {\ell_0^2}r_0$. Then, since
$n={\ell_3}{r_0}+p$ where $\ell_3 \geq \ell_0^2, 0\leq p< r_0$, by
Theorem \ref{thm3} we have $rc_k(G)=2$, here $\ell_0= \lceil \max
\{\sqrt{\frac{k}{2}+1},\frac{k-1}{r}+2 \}\rceil = \lceil \max
\{\sqrt{\frac{k}{2}+1},\frac{k-1}{\lceil \sqrt{2k}\rceil}+2
\}\rceil$=$\lceil \frac{k-1}{\lceil \sqrt{2k}\rceil} \rceil+2$.
Thus we have

\begin{cor}\label{cor1}
For every integer $k \geq 2$, let integer $f(k)={\ell_0^2}r_0$. If
$n \geq f(k)$, then $rc_k(K_n)=2,$ where $\ell_0= \lceil
\frac{k-1}{\lceil \sqrt{2k}\rceil} \rceil+2$ and $r_0=\lceil
\sqrt{2k}\rceil$. \qed
\end{cor}

From Corollary \ref{cor1}, we know that for the integer $f(k)$, we
have that $f(k)\leq ck^{\frac{3}{2}}+C$ where $0< c< 1$ is a
constant and $C=o(k^{\frac{3}{2}})$, and this upper bound is a
considerable improvement of the bound given in Theorem \ref{thm1} of
\cite{Chartrand 2}.

\end{document}